\pdfoutput=1
\documentclass[11pt,letterpaper]{article}
\usepackage[margin=1in]{geometry}
\usepackage[T1]{fontenc}
\usepackage{amsmath,amsthm,amssymb,mathtools}
\usepackage{mathptmx}
\usepackage[protrusion=false]{microtype}
\usepackage{setspace}
\usepackage[authoryear,round]{natbib}
\bibpunct[, ]{(}{)}{,}{a}{}{,}
\usepackage{graphicx,tikz,booktabs,needspace}
\usetikzlibrary{patterns}
\usepackage[hidelinks]{hyperref}
\usepackage{xurl}
\usepackage{bookmark}
\makeatletter
\g@addto@macro\normalsize{%
 \setlength{\abovedisplayskip}{5pt plus 2pt minus 2pt}\setlength{\belowdisplayskip}{5pt plus 2pt minus 2pt}%
 \setlength{\abovedisplayshortskip}{2pt plus 2pt minus 1pt}\setlength{\belowdisplayshortskip}{3pt plus 2pt minus 2pt}}
\renewcommand\section{\@startsection{section}{1}{\z@}{-2ex plus -.5ex minus -.2ex}{1ex plus .2ex}{\normalfont\Large\bfseries\setstretch{1}}}
\renewcommand\subsection{\@startsection{subsection}{2}{\z@}{-1.5ex plus -.4ex minus -.2ex}{.6ex plus .1ex}{\normalfont\large\bfseries\setstretch{1}}}
\renewcommand\paragraph{\@startsection{paragraph}{4}{\z@}{1ex plus .3ex minus .2ex}{-1em}{\normalfont\normalsize\bfseries}}
\g@addto@macro\appendix{\renewcommand{\@seccntformat}[1]{Appendix~\csname the#1\endcsname\quad}}
\makeatother
\allowdisplaybreaks[1]
\newcommand{\pos}[1]{[#1]_+}
\newcommand{\E}{\mathbb E}
\newcommand{\Prob}{\mathbb P}
\newcommand{\dd}{\,d}
\newcommand{\norm}[1]{\lVert#1\rVert}
\newtheoremstyle{statement}{5pt plus 2pt minus 1pt}{5pt plus 2pt minus 1pt}{\itshape}{}{\bfseries}{.}{.5em}{}
\theoremstyle{statement}
\newtheorem{theorem}{Theorem}
\newtheorem{lemma}{Lemma}
\newtheorem{proposition}{Proposition}
\newtheorem{corollary}{Corollary}
\makeatletter
\renewenvironment{proof}[1][\proofname]{\par\pushQED{\qed}\normalfont\topsep4\p@\@plus2\p@\@minus1\p@\relax
 \trivlist\item[\hskip\labelsep\itshape#1\@addpunct{.}]\ignorespaces}{\popQED\endtrivlist\@endpefalse}
\makeatother
\hypersetup{pdftitle={When Speed Changes Where to Build: Optimal Placement of an Optional Fast Connection},pdfauthor={Heng Chen, Xiaofeng Nie},pdfsubject={Optimal terminal placement, optional routing, and the value of speed}}
\begin{document}
\begin{center}
{\Large\bfseries When Speed Changes Where to Build:\\ Optimal Placement of an Optional Fast Connection}\par
\vspace{10pt}
{\setstretch{1.12}%
\begin{minipage}[t]{0.46\linewidth}\centering
{\large Heng Chen}\par\vspace{2pt}
{\small College of Business\\ Kansas State University\\ Manhattan, KS 66506, USA\\
\href{mailto:hengchen@ksu.edu}{\nolinkurl{hengchen@ksu.edu}}\par}
\end{minipage}\hfill
\begin{minipage}[t]{0.50\linewidth}\centering
{\large Xiaofeng Nie}\par\vspace{2pt}
{\small Broadwell College of Business and Economics\\ Fayetteville State University\\ Fayetteville, NC 28301, USA\\
\href{mailto:xnie@uncfsu.edu}{\nolinkurl{xnie@uncfsu.edu}}\par}
\end{minipage}\par}
\vspace{2pt}
\end{center}
\noindent\textbf{Abstract.}
Express roads, limited-stop transit, and air taxis speed up only the trips that choose them, while others stay on the ground network. Where should such a connection's terminals go? Because its users change with the design, expected travel time is nonconvex and nonsmooth in the terminal locations. Previous models have therefore been solved numerically or in restricted settings. We solve a canonical case exactly: uniform trips in a square with grid streets and a straight link. We find that every optimal design is centered and diagonal or axis aligned, and the alignment switches once, near 2.18 times ground speed. The two designs tied there differ sharply in coverage and ridership. A planner also choosing speed, at a smooth cost of speed, never selects the switching speed. If the cost is proportional to speed, above a threshold unit cost the planner builds a slower diagonal link, below it a faster axis-aligned one, and the optimal speed skips a band around the switch. On a rectangle, every optimum is again centered, centered designs have a closed-form travel time, and only two families of links can be optimal, those along the long axis and those tilted toward the long side.

\medskip
\noindent\textbf{Subject classifications:} Transportation: models: location (terminals of an optional express link); facilities/equipment planning: location: continuous (two-terminal placement). \textbf{Area of review:} Transportation. \textbf{Keywords:} continuous hub location; route choice; service coverage.

\section{Introduction}
Many transportation improvements are optional. An express road, a limited-stop transit line, an airport shuttle, or an air-taxi route speeds up the trips that choose to use it, while every other traveler continues on the slower ground network. A traveler must first reach an entry terminal and, after leaving the connection, complete the remaining trip. Where the terminals are placed therefore determines both how much time the connection saves and who uses it. Because the set of users changes with the design, expected travel time is a nonconvex and nonsmooth function of the terminal locations (see the discussion after Theorem~\ref{thm:master}). Previous models have therefore been solved numerically, within parametric families, or for finite trip sets (Section~\ref{sec:lit}).

Urban air mobility is a motivating application. Air-taxi services need accessible takeoff and landing facilities, often called vertiports, and their benefit depends on total passenger travel time rather than aircraft speed alone. The Federal Aviation Administration's (FAA) airspace blueprint expects that, as operations grow, air taxis will fly along corridors linking major airports with city-center vertiports \citep{faa2023}. Research on unmanned aerial vehicle operations has examined routing, uncertainty, and economic, environmental, and social costs \citep{hu2024}, whereas we ask where the terminals should be placed. Also, express road interchanges, limited-stop transit stations, underground loops such as the Las Vegas Convention Center Loop \citep{boringloop}, and airport moving walkways \citep{schwerdfeger2024} pose the same access--travel--egress choice. The technologies differ in construction and operating costs and in user experience, but in each, the faster link is not the fastest \emph{door-to-door} option for every trip.

We study the simplest version of the problem that retains this difficulty: a single bidirectional fast connection with two freely located terminals, travelers who take the fastest door-to-door route, and trips dispersed over a district. Three questions follow. First, where should the terminals go, and how does the answer change with the connection's speed? Second, how much of the district can benefit, what fraction of trips use the link, and what is additional speed worth? Third, how do these answers change with the shape of the district? We answer them in a continuous spatial model in which origins and destinations are independent and uniformly distributed, ground travel follows the Manhattan metric of a street grid, the connection carries travelers faster along the straight line between its terminals, and every traveler may bypass it.

Our main finding is that speed changes where the connection should be built. For uniform demand on a square, we solve the location problem exactly. Every optimal terminal pair is symmetric about the center, the expected travel time of every symmetric pair has a closed form, and only a diagonal or an axis-aligned connection can be optimal. The optimal orientation changes once, from diagonal to axis aligned, at a relative speed of about 2.183, and the change is discontinuous in the terminal locations. It reflects a shift in the source of the connection's value. At speeds near ground speed, an axis-aligned connection saves little time, being as long as the street path it replaces, whereas a diagonal one still saves time by cutting across the grid. As speed grows, speed itself dominates, and the axis-aligned layout, which more trips can use, overtakes the diagonal.

Beyond location, we derive the geographical coverage, the proportion of trips using the connection, and the marginal value of speed in closed form. The two layouts that tie at the switching speed have the same expected travel time but serve very different travelers: about $58\%$ of origins and $12\%$ of trips for the diagonal layout, against $78\%$ and $19\%$ for the axis-aligned layout. The axis-aligned layout carries more trips and gains more from additional speed. The choice of technology and the choice of location therefore depend on each other. A planner who also chooses the speed, at any smooth cost of speed, never chooses the switching speed. With a cost proportional to speed, a single cost level separates the two alignments, and a band of speeds around the switching speed is never chosen.

The analysis extends to a rectangular district. Every optimal terminal pair is again symmetric about the center, and the expected travel time of every symmetric pair again has a closed form. Only two families of connections can be optimal, those along the long axis and those tilted toward the long side.

The contributions of this paper are twofold. First, we show that faster technology can change where an optional connection should be built, and that this change determines which travelers are served. Second, to the best of our knowledge, our paper provides the first complete global characterization of the optimal, freely located terminals of an optional straight Euclidean connection, for independent uniform trips on a square with Manhattan ground travel. The closed forms also provide planners with benchmarks for coverage, ridership, and the value of speed.

\section{Related literature}\label{sec:lit}
Our work is related to three streams of research. The first stream is hub and link location. The classical models of \citet{okelly1986,okelly1987} and the hub-arc models of \citet{cek2005a,cek2005b} address routed-flow consolidation. \citet{campbellokelly2012} and \citet{alumur2021} review this literature. Routing through hubs is usually mandatory in these models. An exception is the nonstrict hubbing of \citet{aykin1995}, in which a hub is used only when this is cheaper. Continuous-demand models instead examine flows dispersed over a region (\citealp{hall1984}; \citealp{campbell1990freight,campbell1990terminals,campbell1992,campbell1993}; \citealp{carlsson2013}). In transit planning, continuum models represent a city as a square with uniform demand and choose among parametric families of networks \citep{daganzo2010}. Our setting combines the continuous spatial demand of these models with a \emph{choice} between a direct trip and travel through an accelerated inter-hub link. In the closely related Manhattan model of \citet{wu2026}, providers must visit a hub before serving customers. Here, no facility visit is required, and the possibility of bypassing the express connection generates the coverage, ridership, and speed-dependent orientation questions. For urban air mobility, vertiport location is mainly formulated as discrete hub-location or network-planning optimization \citep{willey2021,rath2022,wang2022}.

The second stream addresses the geometric placement of highways and other devices that accelerate travel. \citet{abellanas2003} and \citet{baechwa2006} study transportation-induced distance structures, whereas \citet{korman2008} and \citet{ahn2009} consider where to place a high-speed highway in the plane, including under the city metric. The closest of these models is the moving walkway of \citet{cardinal2008}: two points joined by a fast link that a traveler uses only when it helps, located to minimize the maximum travel time between the points of a given finite set. Related models locate a high-speed line jointly with a service facility (\citealp{espejo2011}; \citealp{diazbanez2013,diazbanez2014,diazbanez2016}). In contrast, we minimize the expected \emph{door-to-door} time over continuously distributed origin--destination pairs, and we focus on the globally optimal endpoints of a single optional link, whose Euclidean length, and hence travel time, is determined endogenously by the placement decision.

The third stream treats rapid transit and accelerated passenger networks. \citet{gendreau1995} and \citet{bruno2002} examine the location of rapid-transit lines, \citet{laporte2005} construct a single alignment that maximizes the origin--destination demand it covers, \citet{marinjaramillo2009} design lines and stations jointly, and \citet{laporte2011} review rapid-transit network planning. Speed-up subnetworks, accelerating moving walkways, and airport-terminal walkways provide related settings in which access, network design, and service speed interact \citep{schmidt2014,scarinci2017,boysen2021,schwerdfeger2024}, and \citet{korner2014} locate entry and exit points on a given faster network to maximize the number of covered origin--destination pairs. In these models the alignments are chosen among discrete candidates or on a given network, the objective is covered demand, ridership, or cost, and speed enters as an input. We instead locate the terminals freely in a continuous region, minimize expected door-to-door time, and show how speed alone changes the optimal alignment.

\section{Model and centering}\label{sec:centering}
Let the origin $X=(X_1,X_2)$ and destination $Y=(Y_1,Y_2)$ be independent uniform points of $\Omega=[0,1]^2$, and write $d(X,Y)=\norm{X-Y}_1$ for ground distance. Ground speed is one. Two distinct terminals $P,Q\in\Omega$ are joined by a straight bidirectional link of length $\ell=\norm{Q-P}_2$ and relative speed $v>1$. Put $q=1/v$, the inverse relative speed, and $L=q\ell$, the link travel time. We allow $q=0$ as the instantaneous-link limit. Trips take their fastest routes, and the planner minimizes expected travel time:
\begin{equation}\label{eq:model}
 F=\E\min\{d(X,Y),\ d(X,P)+L+d(Q,Y),\ d(X,Q)+L+d(P,Y)\}.
\end{equation}
The three terms are the travel times of a trip that stays on the ground, rides the link from $P$ to $Q$, or rides it from $Q$ to $P$. Continuity and compactness ensure that a minimizing pair exists. Define the forward gain $g_{PQ}=d(X,Y)-d(X,P)-L-d(Q,Y)$, the time saved by riding from $P$ to $Q$, and define $g_{QP}$ by exchanging $P$ and $Q$. Bounding $d(X,Y)$ by $d(X,Q)+d(Q,Y)$ in $g_{PQ}$, and by $d(X,P)+d(P,Y)$ in $g_{QP}$, gives
\begin{equation}\label{eq:access}
 g_{PQ}\le d(X,Q)-d(X,P)-L,\qquad
 g_{QP}\le d(X,P)-d(X,Q)-L.
\end{equation}
The two bounds sum to $-2L\le0$, so the two gains are never both positive, and a trip saves $[g_{PQ}]_++[g_{QP}]_+$, where $[z]_+=\max(z,0)$. Exchanging $X$ and $Y$ turns $g_{QP}$ into $g_{PQ}$, so the two terms have equal expectations, and $\E d(X,Y)=2/3$ yields
\begin{equation}\label{eq:saving}
 F=\frac23-2\E[g_{PQ}]_+.
\end{equation}
By \eqref{eq:access} and $d(X,Q)\le d(X,P)+d(P,Q)$, no gain exceeds $d(P,Q)-L$, which is the gain of the trip from $P$ to $Q$. This gain is positive, because $d(P,Q)=\norm{Q-P}_1\ge\norm{Q-P}_2>L$, and trips that start near $P$ and end near $Q$ gain as well. Hence $F<2/3$ for every pair of terminals. Ground distance is a sum over the two coordinates, so we first study the gain in one coordinate.

\begin{lemma}[Centering a coordinate gain]\label{lem:window}
Fix $0\le\delta\le1$ and $|t|\le(1-\delta)/2$. Two terminals lie in $[0,1]$ at $(1-\delta)/2+t$ and $(1+\delta)/2+t$: they are $\delta$ apart and shifted by $t$ from the centered position. For $u$ and $w$ independent and uniform on $[0,1]$, $Z=|w-u|-{|u-(1-\delta)/2-t|}-{|w-(1+\delta)/2-t|}$ is the coordinate gain of the trip from $u$ to $w$: the ground distance it saves in this coordinate by riding from the first terminal to the second, before link time is charged. The probability $\Lambda_z(t)=\Prob(Z>z)$ that the gain exceeds $z$ is
\begin{equation}\label{eq:window}
 \Lambda_z(t)=
 \begin{cases}
 0,&z\ge\delta,\\
 (1-z)^2/4-(\delta-z)^2/8-t^2,&-\delta\le z<\delta,\\
 1-[(1+z)/2-t]_+^2-[(1+z)/2+t]_+^2,&z<-\delta.
 \end{cases}
\end{equation}
Hence $\Lambda_z(t)\le\Lambda_z(0)$ for every $z$, with strict inequality whenever $t\ne0$ and $-1<z<\delta$.
\end{lemma}
\begin{figure}[t]
\setstretch{1.5}
\centering
\begin{tikzpicture}[x=4.4cm,y=4.4cm,font=\footnotesize,>=stealth]
\def\fga{0.35}\def\fgb{0.75}
\node[anchor=south west,inner sep=0pt] at (0,1.2) {\textbf{(a) Ranges of the gain $Z$}};
\fill[black!25] (0,\fga)--(\fga,\fga)--(\fgb,\fgb)--(\fgb,1)--(0,1)--cycle;
\fill[black!50] (0,\fgb) rectangle (\fga,1);
\fill[black!9] (0,0) rectangle (\fga,\fga) (\fgb,\fgb) rectangle (1,1);
\fill[pattern=north east lines,pattern color=black!60] (\fga,0)--(1,0)--(1,\fgb)--(\fgb,\fgb)--(\fga,\fga)--cycle;
\draw[very thin,densely dashed] (\fga,0)--(\fga,\fga) (\fga,\fgb)--(\fga,1) (\fgb,\fgb)--(\fgb,1) (0,\fga)--(\fga,\fga) (0,\fgb)--(\fga,\fgb) (\fgb,\fgb)--(1,\fgb);
\draw[very thin] (\fga,\fga)--(\fgb,\fgb);
\draw (0,0) rectangle (1,1);
\node[fill=white,inner sep=1.5pt] at (.175,.875) {$Z=\delta$};
\node at (.23,.55) {$-\delta<Z<\delta$};
\node[fill=white,inner sep=1.5pt] at (.72,.27) {$Z=-\delta$};
\node at (.175,.175) {$Z<-\delta$};
\node[font=\scriptsize,inner sep=0pt] at (.875,.875) {$Z<-\delta$};
\foreach \p/\l in {0/0,\fga/\xi,\fgb/\xi+\delta,1/1} {\draw (\p,0)--(\p,-.02) node[below,inner sep=2pt]{$\l$};}
\foreach \p/\l in {\fga/\xi,\fgb/\xi+\delta,1/1} {\draw (0,\p)--(-.02,\p) node[left,inner sep=2pt]{$\l$};}
\node[below] at (.5,-.085) {origin $u$};
\node[rotate=90] at (-.27,.5) {destination $w$};
\begin{scope}[shift={(1.72,0)}]
\node[anchor=south west,inner sep=0pt] at (0,1.2) {\textbf{(b) Trips with $Z>z$} ($-\delta\le z<\delta$)};
\fill[black!25] (0,.5)--(\fga,.5)--(.6,\fgb)--(.6,1)--(0,1)--cycle;
\draw[very thin,densely dashed] (\fga,0)--(\fga,1) (\fgb,0)--(\fgb,1) (0,\fga)--(1,\fga) (0,\fgb)--(1,\fgb);
\draw[line width=.9pt] (0,.5)--(\fga,.5)--(.6,\fgb)--(.6,1);
\draw[line width=.6pt,dash pattern=on 2pt off 1.5pt] (\fga,.5)--(.6,.5)--(.6,\fgb);
\draw (0,0) rectangle (1,1);
\node[fill=white,inner sep=1.5pt] at (.175,.875) {$Z>z$};
\draw[very thin] (.6,1)--(.6,1.1) (0,1)--(0,1.1);
\draw[<->] (0,1.07)--(.6,1.07) node[midway,above,inner sep=2pt]{$\xi+(\delta-z)/2$};
\draw[very thin] (.6,.5)--(1.1,.5) (1,1)--(1.1,1);
\draw[<->] (1.07,.5)--(1.07,1) node[midway,right,inner sep=3pt]{$1-\xi-(\delta+z)/2$};
\node[below,fill=white,inner sep=1.5pt] at (.475,.485) {$(\delta-z)/2$};
\node[right,fill=white,inner sep=1.5pt] at (.615,.64) {$(\delta-z)/2$};
\foreach \p/\l in {0/0,\fga/\xi,\fgb/\xi+\delta,1/1} {\draw (\p,0)--(\p,-.02) node[below,inner sep=2pt]{$\l$};}
\foreach \p/\l in {\fga/\xi,.5/{\xi+(\delta+z)/2},\fgb/\xi+\delta,1/1} {\draw (0,\p)--(-.02,\p) node[left,inner sep=2pt]{$\l$};}
\node[below] at (.5,-.085) {origin $u$};
\end{scope}
\end{tikzpicture}
\caption{The coordinate gain $Z$ of Lemma~\ref{lem:window}, with ${\xi=(1-\delta)/2+t}$, drawn for $\delta=0.4$ and $t=0.05$. (a) Ranges of $Z$. (b) The trips with $Z>z$ when $-\delta\le z<\delta$, drawn for $z=-0.1$: a rectangle minus the dashed triangle. Each point $(u,w)$ is a trip from $u$ to $w$ in one coordinate, not a location in $\Omega$.}\label{fig:gain}
\end{figure}
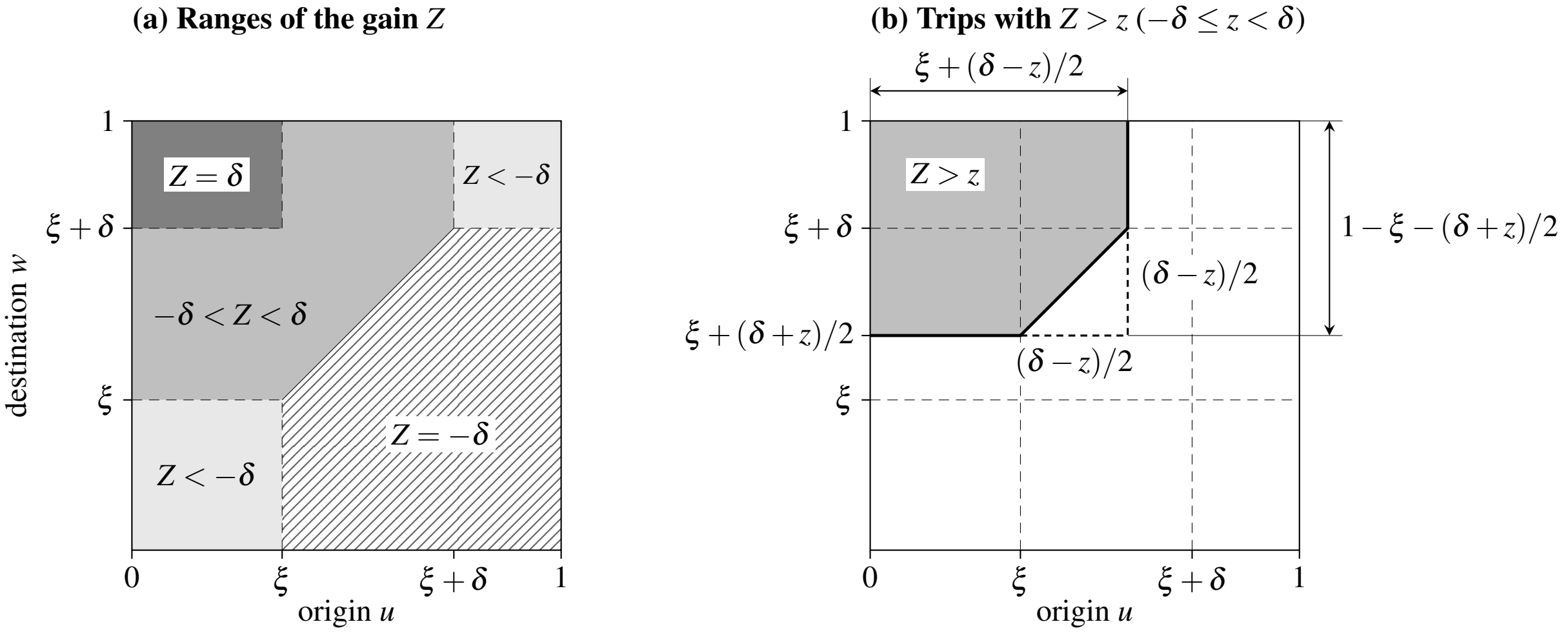
\begin{proof}
Write $\xi=(1-\delta)/2+t$ for the first terminal, so that the second is $\xi+\delta$. After rearranging terms, the gain of a forward trip ($u\le w$) is $Z=\delta-2[u-\xi]_+-2[\xi+\delta-w]_+\le\delta$, whereas that of a backward trip ($u>w$) is $Z=-\delta-2[\xi-u]_+-2[w-\xi-\delta]_+\le-\delta$. Figure~\ref{fig:gain}(a) partitions the square of trips $(u,w)$ according to the resulting ranges of the gain. In particular, no gain exceeds $\delta$, so $\Lambda_z(t)=0$ for $z\ge\delta$.

For $-\delta\le z<\delta$, the trips with $Z>z$ are all forward trips. They fill the rectangle ${u<\xi+(\delta-z)/2}$, ${w>\xi+(\delta+z)/2}$, minus the triangle where ${w\le u+(\delta+z)/2}$ (Figure~\ref{fig:gain}(b)). The rectangle has sides ${\xi+(\delta-z)/2}$ and ${1-\xi-(\delta+z)/2}$, and the triangle has legs $(\delta-z)/2$, so $\Lambda_z(t)=\{\xi+(\delta-z)/2\}\{1-\xi-(\delta+z)/2\}-(\delta-z)^2/8=\{(1-z)/2+t\}\{(1-z)/2-t\}-(\delta-z)^2/8=(1-z)^2/4-t^2-(\delta-z)^2/8$.

For $z<-\delta$, the trips with $Z\le z$ have $u$ and $w$ both below $\xi$, where $Z=-\delta-2(\xi-\max\{u,w\})$, or both above $\xi+\delta$, where $Z=-\delta-2(\min\{u,w\}-\xi-\delta)$. Hence $Z\le z$ exactly when $\max\{u,w\}\le\xi+(z+\delta)/2=(1+z)/2+t$ or $\min\{u,w\}\ge\xi+\delta-(z+\delta)/2=1-\{(1+z)/2-t\}$. These corner squares have sides $[(1+z)/2+t]_+$ and $[(1+z)/2-t]_+$, and subtracting their areas from one gives the last expression.

Finally, for $-\delta\le z<\delta$ the middle case shows $\Lambda_z(t)=\Lambda_z(0)-t^2$. The last expression is even in $t$. For $t>0$ its derivative $2[(1+z)/2-t]_+-2[(1+z)/2+t]_+$ is at most zero, and negative when $z>-1$.
\end{proof}

\begin{theorem}[Every optimal pair is centered]\label{thm:center}
For fixed displacement $\Delta=Q-P$, centering uniquely minimizes $F$ over feasible translations. Every global optimum thus has $P=(1/2,1/2)-\Delta/2$ and $Q=(1/2,1/2)+\Delta/2$.
\end{theorem}
\begin{proof}
By reflection, take $\Delta=(\delta_1,\delta_2)\ge0$, and let $(t_1,t_2)$ be the translation of the pair from the centered position, so that coordinate $k$ of $P$ is $(1-\delta_k)/2+t_k$ and that of $Q$ is $(1+\delta_k)/2+t_k$. Ground distance is a sum over the two coordinates, so $g_{PQ}=Z_1+Z_2-L$, where $Z_k$ is the gain of Lemma~\ref{lem:window} for coordinate $k$ of the trip: the lemma with $u=X_k$, $w=Y_k$, $\delta=\delta_k$, and $t=t_k$. It applies because the four coordinates $X_1,X_2,Y_1,Y_2$ are independent and uniform on $[0,1]$. For the same reason, $Z_1$ and $Z_2$ are independent. Given $Z_2=z_2$, the expected saving is $\E[Z_1+z_2-L]_+=\int_{L-z_2}^{\infty}\Prob(Z_1>z)\dd z$, and by Lemma~\ref{lem:window} the integrand $\Prob(Z_1>z)=\Lambda_z(t_1)$ is largest at $t_1=0$ for every $z$. Hence, for every $t_2$, the expected saving $\E[g_{PQ}]_+$ is at least as large at $t_1=0$ as at any other $t_1$, because this holds given each value $z_2$ and hence on average over $Z_2$. Likewise, for every $t_1$ it is at least as large at $t_2=0$ as at any other $t_2$. By \eqref{eq:saving}, $F$ is therefore smallest when $t_1=t_2=0$.

For uniqueness, take a pair with $t_1\ne0$, relabeling the coordinates if necessary. Moving it to $t_2=0$ does not lower the expected saving. Moving it then to $t_1=0$ raises $\Lambda_z(t_1)$ strictly for $-1<z<\delta_1$ by Lemma~\ref{lem:window}, and hence $\int_{L-z_2}^{\infty}\Lambda_z(t_1)\dd z$ strictly whenever $L-z_2<\delta_1$. This event, $Z_2>L-\delta_1$, has positive probability, because $L-\delta_1<d(P,Q)-\delta_1=\delta_2$ and $\Prob(Z_2>z)=\Lambda_z(0)>0$ for every $z<\delta_2$ by \eqref{eq:window} with $\delta=\delta_2$. The centered pair therefore saves strictly more than every translated pair, so centering is the unique minimizer.
\end{proof}
Theorem~\ref{thm:center} holds at every speed. By symmetry, we may also exchange the two coordinates, so from now on we let ${0\le\delta_2\le\delta_1\le1}$. Thus $\delta_1$ is the larger coordinate displacement and $\delta_2$ the smaller.

\section{An exact objective for centered pairs}\label{sec:objective}
By the proof of Theorem~\ref{thm:center} and \eqref{eq:saving}, ${F=2/3-2\E[Z_1+Z_2-L]_+}$, where $Z_1$ and $Z_2$ are the independent gains of a trip in the two coordinates. This section gives $F$ in closed form for centered pairs, in three steps.

\paragraph{Step 1: the distribution of a coordinate gain.}
Setting $t=0$ in \eqref{eq:window} gives the distribution of the coordinate gain $Z$ of centered terminals separated by $\delta$. Where $\Lambda_z(0)=\Prob(Z>z)$ is smooth in $z$, its negative derivative is the density $f_Z$, called the lower density below $-\delta$ and the central density above it. The two downward jumps of $\Lambda_z(0)$ are point masses:
\begin{equation}\label{eq:law}
\begin{aligned}
 f_Z(z)&=1+z\quad(-1<z<-\delta),&\qquad \Prob(Z=-\delta)&=\tfrac12-(1-\delta)^2/4,\\
 f_Z(z)&=\frac{1-z}{2}-\frac{\delta-z}{4}=\frac{1-\delta}{2}+\frac{\delta-z}{4}\quad(-\delta<z<\delta),&\qquad \Prob(Z=\delta)&=(1-\delta)^2/4.
\end{aligned}
\end{equation}
At $\delta$ the tail probability falls from $(1-\delta)^2/4$ to zero. At $-\delta$ it falls from $1-(1-\delta)^2/2$ (last case of \eqref{eq:window}) to $(1+\delta)^2/4-\delta^2/2$ (middle case), which equals $\frac12-(1-\delta)^2/4$ because $(1+\delta)^2+(1-\delta)^2=2+2\delta^2$. The jump is therefore $\frac12-(1-\delta)^2/4$ as well. In Figure~\ref{fig:gain}(a), the mass at $\delta$ consists of the forward trips that cover the terminal interval and save the full $\delta$, and the mass at $-\delta$ of the backward trips that meet it.

\paragraph{Step 2: regular and irregular placements.}
Each coordinate gain $Z_k$ has the distribution \eqref{eq:law} with $\delta=\delta_k$. In coordinate $k$ the centered terminals lie at $(1-\delta_k)/2$ and $(1+\delta_k)/2$, so each is $(1-\delta_k)/2$ away from the nearest side of the square. Write $A=(1-\delta_1)/2$ for this gap in the first coordinate and $b=(1-\delta_2)/2$ for the gap in the second. Put $D=(\delta_1+\delta_2-L)/2$ and $E=[\delta_1-\delta_2-L]_+/2$. A coordinate gain is at most the displacement, $Z_k\le\delta_k$, with equality for the trip from $P$ to $Q$. No trip therefore gains more than that trip, whose gain is $\delta_1+\delta_2-L=2D$. Because $Z_1\le\delta_1$, a trip's gain $Z_1+Z_2-L$ can be positive only if $Z_2$ exceeds the threshold $L-\delta_1$. Two cases arise, depending on whether this threshold is below $-\delta_2$ or not (Figure~\ref{fig:integration}).

A placement is \emph{regular} if $L\ge \delta_1-\delta_2$, so that the threshold $L-\delta_1$ is at or above $-\delta_2$. Then a trip that loses at least $\delta_2$ in the second coordinate has $Z_2\le-\delta_2\le L-\delta_1$ (Figure~\ref{fig:integration}(a)), so its gain $Z_1+Z_2-L$ is at most zero. Every diagonal placement, where $\delta_1=\delta_2$, is regular. A placement is \emph{irregular} if $L<\delta_1-\delta_2$, so that the threshold is below $-\delta_2$. Then the gain in the first coordinate can offset both the link time and a loss of at least $\delta_2$ in the second, so trips with $Z_2\le-\delta_2$ can benefit too, by up to $\delta_1-\delta_2-L=2E$. Every axis placement, where $\delta_2=0$, is irregular. In practice, when the link time is shorter than the imbalance $\delta_1-\delta_2$ between its two displacements, that is, for a fast link close enough to an axis, travelers can benefit even when their trip moves against the link in the coordinate of the smaller displacement. This changes which trips the link can serve, not where its terminals may be placed.

\paragraph{Step 3: the cost in closed form.}
With $s$ the length of the side carrying $\delta_2$ ($s\ne1$ only for rectangles), define
\begin{align}
 R(A,b,D)&=\frac{D^5}{60}+\frac{(A+b)D^4}{6}
 +\frac{(A^2+4Ab+b^2)D^3}{3}
 +2Ab(A+b)D^2+2A^2b^2D,\label{eq:R}\\
 C_{A,b,s}(E)&=\Bigl(\frac{s^2}{2}-b^2\Bigr)\Bigl(2A^2E+2AE^2+\frac{E^3}{3}\Bigr)
 +\Bigl(2b-\frac s2\Bigr)\Bigl(2A^2E^2+\frac{4AE^3}{3}+\frac{E^4}{6}\Bigr)\notag\\
 &\quad-\frac{5A^2E^3}{3}-\frac{5AE^4}{6}-\frac{E^5}{12}.\label{eq:C}
\end{align}

\begin{theorem}[Exact centered objective]\label{thm:master}
For every centered terminal placement on the unit square, $G:=\E[g_{PQ}]_+=R(A,b,D)+C_{A,b,1}(E)$ and $F=2/3-2G$. In the regular regime $E=0$, so that $G=R(A,b,D)$.
\end{theorem}
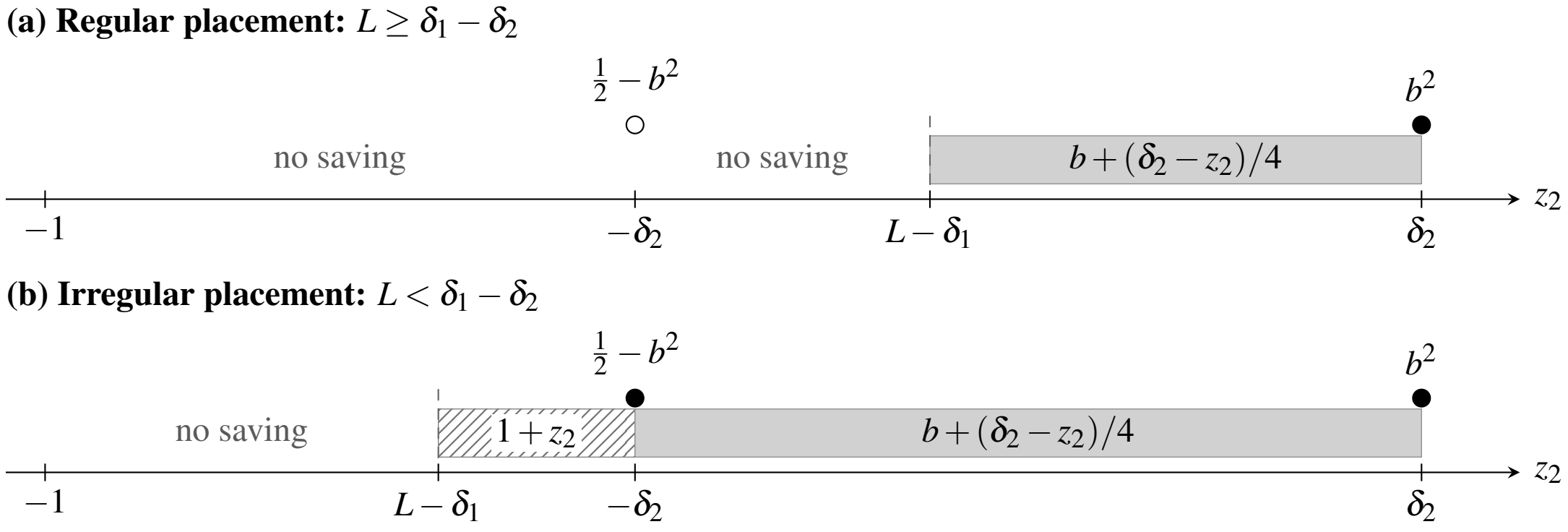
\begin{figure}[t]
\setstretch{1.5}
\centering
\begin{tikzpicture}[x=9cm,y=1cm,font=\footnotesize,>=stealth]
\tikzset{bar/.style={fill=black!18,draw=black!45,line width=.25pt},mass/.style={circle,draw=black,fill=black,inner sep=0pt,minimum size=4.5pt},nomass/.style={circle,draw=black,fill=white,inner sep=0pt,minimum size=4.5pt},cut/.style={densely dashed,black!60}}
\begin{scope}[yshift=2.5cm]
\node[anchor=west,inner sep=0pt] at (-1.04,1.62) {\textbf{(a) Regular placement:} $L\ge\delta_1-\delta_2$};
\draw[->] (-1.04,0)--(0.5,0) node[right] {$z_2$};
\draw[bar] (-0.1,.14) rectangle (0.4,.58);
\node at (0.15,.36) {$b+(\delta_2-z_2)/4$};
\node[black!65] at (-0.7,.36) {no saving};
\node[black!65] at (-0.25,.36) {no saving};
\foreach \p/\l in {-1/{-1},-0.4/{-\delta_2},-0.1/{L-\delta_1},0.4/{\delta_2}} {\draw (\p,.08)--(\p,-.08) node[below,inner sep=2pt] {$\l$};}
\draw[cut] (-0.1,.13)--(-0.1,.78);
\node[nomass] at (-0.4,.68) {};\node[above] at (-0.4,.76) {$\tfrac12-b^2$};
\node[mass] at (0.4,.68) {};\node[above] at (0.4,.76) {$b^2$};
\end{scope}
\begin{scope}
\node[anchor=west,inner sep=0pt] at (-1.04,1.62) {\textbf{(b) Irregular placement:} $L<\delta_1-\delta_2$};
\draw[->] (-1.04,0)--(0.5,0) node[right] {$z_2$};
\draw[pattern=north east lines,pattern color=black!55,draw=black!50,line width=.25pt] (-0.6,.14) rectangle (-0.4,.58);
\draw[bar] (-0.4,.14) rectangle (0.4,.58);
\node[fill=white,inner sep=1pt] at (-0.5,.36) {$1+z_2$};
\node at (0,.36) {$b+(\delta_2-z_2)/4$};
\node[black!65] at (-0.8,.36) {no saving};
\foreach \p/\l in {-1/{-1},-0.6/{L-\delta_1},-0.4/{-\delta_2},0.4/{\delta_2}} {\draw (\p,.08)--(\p,-.08) node[below,inner sep=2pt] {$\l$};}
\draw[cut] (-0.6,.13)--(-0.6,.78);
\node[mass] at (-0.4,.68) {};\node[above] at (-0.4,.76) {$\tfrac12-b^2$};
\node[mass] at (0.4,.68) {};\node[above] at (0.4,.76) {$b^2$};
\end{scope}
\end{tikzpicture}
\caption{The parts of the distribution of $Z_2$ integrated in the proof of Theorem~\ref{thm:master}, for $\delta_1=0.8$, $\delta_2=0.4$, and $L=0.7$ in (a) or $L=0.2$ in (b). Bars are the central density (gray) and the lower density (hatched) above the threshold $L-\delta_1$. A filled dot is a point mass that contributes, an open dot one that does not.}\label{fig:integration}
\end{figure}
\begin{proof}
We integrate first over $Z_1$, given $Z_2=z_2$. Since $z_2\le\delta_2\le\delta_1$, no value $Z_1\le-\delta_1$ gives a positive gain. By \eqref{eq:law}, the rest of the distribution of $Z_1$ is the mass $A^2$ at $\delta_1$ and the central density $A+(\delta_1-z_1)/4$ on $(-\delta_1,\delta_1)$. Let $\Phi(z_2)=\E[Z_1+z_2-L]_+$ be the expected saving given $Z_2=z_2$. It is zero for $z_2\le L-\delta_1$. For $z_2>L-\delta_1$, the gain $z_1+z_2-L$ is positive exactly when $z_1>L-z_2$, and ${L-z_2\ge-\delta_2\ge-\delta_1}$ lies in the range of the central density. The mass and the part of the central density above $L-z_2$ therefore give
\[
 \mbox{\footnotesize$\displaystyle\Phi(z_2)=A^2(\delta_1+z_2-L)+\int_{L-z_2}^{\delta_1}(z_1+z_2-L)\Bigl(A+\frac{\delta_1-z_1}{4}\Bigr)\dd z_1
 =A^2(\delta_1+z_2-L)+\frac A2(\delta_1+z_2-L)^2+\frac{(\delta_1+z_2-L)^3}{24}.$}
\]

We next integrate over $Z_2$, because $G=\E\,\Phi(Z_2)$ by independence. By \eqref{eq:law}, the distribution of $Z_2$ has the masses $b^2$ at $\delta_2$ and $\frac12-b^2$ at $-\delta_2$, the central density $b+(\delta_2-z_2)/4$ on $(-\delta_2,\delta_2)$, and the lower density $1+z_2$ on $(-1,-\delta_2)$. Figure~\ref{fig:integration} shows which parts contribute. In both regimes the threshold is $z_2=L-\delta_1$.

\emph{Regular placements} ($L\ge\delta_1-\delta_2$). Here $L-\delta_1\ge-\delta_2$, so neither the lower density nor the mass at $-\delta_2$ contributes (Figure~\ref{fig:integration}(a)). Only the mass at $\delta_2$ and the central density above the threshold remain:
\[
 G=b^2\Phi(\delta_2)+\int_{L-\delta_1}^{\delta_2}\Bigl(b+\frac{\delta_2-z_2}{4}\Bigr)\Phi(z_2)\dd z_2=R(A,b,D).
\]
The last equality follows by inserting $\Phi$ and integrating, with $\delta_1+\delta_2-L=2D$.

\emph{Irregular placements} ($L<\delta_1-\delta_2$). Now $L-\delta_1<-\delta_2$, so the mass at $-\delta_2$ and part of the lower density contribute as well (Figure~\ref{fig:integration}(b)). The four contributions add up to
\begin{equation}\label{eq:master-integral}
 G=b^2\Phi(\delta_2)+\int_{-\delta_2}^{\delta_2}\Bigl(b+\frac{\delta_2-z_2}{4}\Bigr)\Phi(z_2)\dd z_2+\Bigl(\frac12-b^2\Bigr)\Phi(-\delta_2)+\int_{L-\delta_1}^{-\delta_2}(1+z_2)\Phi(z_2)\dd z_2.
\end{equation}
To connect this expression with $R$, observe that the last equality of the regular case does not use $L\ge\delta_1-\delta_2$. Hence $R(A,b,D)$ is the value obtained by extending the central density down to the threshold $L-\delta_1$. Below $-\delta_2$, however, $Z_2$ has the lower density $1+z_2$, and at $-\delta_2$ it has the mass $\frac12-b^2$. The correction is therefore
\[
 G-R(A,b,D)=\Bigl(\frac12-b^2\Bigr)\Phi(-\delta_2)+\int_{L-\delta_1}^{-\delta_2}\Bigl\{1+z_2-\Bigl(b+\frac{\delta_2-z_2}{4}\Bigr)\Bigr\}\Phi(z_2)\dd z_2.
\]
The interval from $L-\delta_1$ to $-\delta_2$ has length $\delta_1-\delta_2-L=2E$ and lies in the range of the lower density because $\delta_1-L\le1$. Inserting $\Phi$ and integrating gives $C_{A,b,1}(E)$. Hence $G=R(A,b,D)+C_{A,b,1}(E)$ in both regimes, with $E=0$ for regular placements. If $\delta_2=0$, the two masses of $Z_2$ coincide and add. Finally, \eqref{eq:saving} gives $F$.
\end{proof}
The correction $C_{A,b,1}(E)$ is present exactly when $L<\delta_1-\delta_2$. With $L=q\sqrt{\delta_1^2+\delta_2^2}$, this reads ${(\delta_1-\delta_2)/\sqrt{\delta_1^2+\delta_2^2}>q}$, that is, for links close enough to an axis. The trips with $(Z_1,Z_2)=(\delta_1,-\delta_2)$, in the regions $Z=\delta$ and $Z=-\delta$ of Figure~\ref{fig:gain}(a), all have the gain $\delta_1-\delta_2-L$, and for $A>0$ their probability, $A^2(\frac12-b^2)$, is positive. They all start to ride at once when $L$ falls below $\delta_1-\delta_2$, and from there on the cost falls faster as $L$ decreases, so it has a kink at this value. Nor is the cost convex in the terminal locations. Let the first terminal move from $P$ to $Q$ while the second moves from $Q$ to $P$. The cost is the same at both ends, below $2/3$, but halfway the two terminals meet and the link saves nothing, so the cost is $2/3$.

\section{Globally optimal placement on the square}\label{sec:square}
It remains to choose $\delta_1$ and $\delta_2$. A centered pair is \emph{diagonal} if $\delta_2=\delta_1$, an \emph{axis} pair if $\delta_2=0$, and \emph{tilted} if $0<\delta_2<\delta_1$. Theorem~\ref{thm:orient} rules out the tilted case.

\begin{theorem}[Only diagonal and axis layouts can be optimal]\label{thm:orient}
Up to the symmetries of the square and the exchange of $P$ and $Q$, every global optimum is a centered diagonal pair or a centered axis pair.
\end{theorem}
\begin{proof}
\emph{The saving of a regular placement increases with $\alpha=2Ab$ and with $D$:} Here $G=R(A,b,D)$ and $D>0$. In \eqref{eq:R}, $A^2+4Ab+b^2=(A+b)^2+\alpha$ and $2A^2b^2=\alpha^2/2$, so $R$ depends on $A$ and $b$ only through $A+b$ and $\alpha$. It increases with $D$, and at fixed $A+b$ and $D$ its derivative in $\alpha$ is ${R_\alpha=D^3/3+(A+b)D^2+\alpha D>0}$.

\emph{A regular optimum is diagonal:} If $\delta_1>\delta_2$, reduce $\delta_1-\delta_2$ at fixed $\delta_1+\delta_2$. Then $A+b=1-(\delta_1+\delta_2)/2$ is unchanged, and $\alpha=(1-\delta_1)(1-\delta_2)/2$ increases strictly, because $\delta_1\delta_2=[(\delta_1+\delta_2)^2-(\delta_1-\delta_2)^2]/4$ does. Moreover, $\delta_1^2+\delta_2^2=[(\delta_1+\delta_2)^2+(\delta_1-\delta_2)^2]/2$ decreases, so the link shortens and $D$ does not decrease. The placement stays regular, that is, $(\delta_1-\delta_2)^2/(\delta_1^2+\delta_2^2)\le q^2$, because this ratio equals $2(\delta_1-\delta_2)^2/[(\delta_1+\delta_2)^2+(\delta_1-\delta_2)^2]$ and so decreases with $\delta_1-\delta_2$. By the first step, the saving increases strictly until $\delta_1=\delta_2$.

\emph{An irregular optimum is an axis pair:} The irregular placements form the set $0\le \delta_2<\delta_1-L$, $\delta_1\le1$. It does not include the interface $L=\delta_1-\delta_2$, whose placements are regular. An irregular optimum therefore (a) is an axis pair ($\delta_2=0$), (b) has $\delta_1=1$ and $\delta_2>0$, or (c) lies in the interior of the set ($\delta_2>0$ and $\delta_1<1$). In case (c) it is a stationary point of $G$, because $G$ is differentiable in the interior. Appendix~\ref{app:exclude} excludes (b) and (c). It treats a $\rho\times1$ rectangle and is read here with $\rho=1$ and $(x,y)=(\delta_1,\delta_2)$. Its Part~2 shows that no stationary point as in (c) exists, and its Part~3 (face $x=\rho$) that every placement as in (b) can be improved.
\end{proof}

\paragraph{Diagonal pairs.}
Write a diagonal pair as ${P=(a,a)}$, ${Q=(1-a,1-a)}$, where the offset $a$ is the distance of each terminal from the nearest sides of the square. Both displacements equal $1-2a$ and $L=\sqrt2\,q(1-2a)$. Hence $A=b=a$ and $D=m(1-2a)=:D_d$, where $m=1-q/\sqrt2$ is the share of the bypassed ground distance that the link saves. This pair is regular, so its total saving is
\begin{equation}\label{eq:Sd}
 S_d=2R(a,a,D_d)=\frac{D_d^5}{30}+\frac{2aD_d^4}{3}
 +4a^2D_d^3+8a^3D_d^2+4a^4D_d.
\end{equation}

\begin{proposition}[Unique diagonal optimum]\label{prop:roots}
For $q\in[0,1)$, the diagonal saving \eqref{eq:Sd} is maximized by a unique offset $a\in(0,1/2)$. It satisfies $2a/(1-2a)=\lambda_d$, where $\lambda_d$ is the unique positive root, a simple one, of
\begin{equation}\label{eq:diagroot}
 \Xi_d(\lambda)=3\lambda^4+12(2m-1)\lambda^3+36m(m-1)\lambda^2+8m^2(2m-3)\lambda+2m^3(m-2).
\end{equation}
\end{proposition}
\begin{proof}
Let $\lambda=2a/(1-2a)$. Since $dD_d/da=-2m$, differentiating \eqref{eq:Sd} and inserting $a=\lambda(1-2a)/2$ and $D_d=m(1-2a)$, which makes every term a multiple of $(1-2a)^4$, gives $dS_d/da=-(m/6)(1-2a)^4\,\Xi_d(\lambda)$. The coefficients of $\Xi_d$ have the signs $(+,\pm,-,-,-)$, where the second can have either sign and the third vanishes at $q=0$. In every case they change sign exactly once. By Descartes' rule of signs, $\Xi_d$ has exactly one positive root, counting multiplicity, and is negative below it and positive above it. Since $\lambda$ increases from zero to infinity as $a$ increases from $0$ to $1/2$, the saving first increases and then decreases.
\end{proof}

\paragraph{Axis pairs.}
Write an axis pair as ${P=(p,1/2)}$, ${Q=(1-p,1/2)}$, with offset $p$. Then $\delta_1=1-2p$, $\delta_2=0$, and $L=q(1-2p)$, so $A=p$, $b=1/2$, and $D=E=n(1-2p)/2=:D_a$, where $n=1-q$ is again the share of the bypassed ground distance that the link saves. For $q<1$ the pair is irregular, so its total saving is ${2[R(p,\frac12,D_a)+C_{p,1/2,1}(D_a)]}$, that is,
\begin{equation}\label{eq:Sa}
 S_a=2p^2D_a+2p(1+2p)D_a^2+\frac{1+8p-8p^2}{3}D_a^3
 +\frac{1-4p}{3}D_a^4-\frac{2}{15}D_a^5.
\end{equation}

\begin{proposition}[Unique axis optimum]\label{prop:axisroots}
For $q\in[0,1)$, the axis saving \eqref{eq:Sa} is maximized by a unique offset $p\in(0,1/2)$. It satisfies $2p/(1-2p)=\lambda_a$, where $\lambda_a$ is the unique positive root, a simple one, of
\begin{equation}\label{eq:axisroot}
 \Xi_a(\lambda)=12\lambda^4+48n\lambda^3+18(n^2+2n-2)\lambda^2
 -4(n^3-9n^2+6n+6)\lambda-n(n^3-6n^2+2n+12).
\end{equation}
\end{proposition}
\begin{proof}
Let $\lambda=2p/(1-2p)$. Since $dD_a/dp=-n$, differentiating \eqref{eq:Sa} and inserting $p=\lambda(1-2p)/2$ and $D_a=n(1-2p)/2$ gives $dS_a/dp=-(n/24)(1-2p)^4\,\Xi_a(\lambda)$, once the terms of lower degree are multiplied by powers of $(1+\lambda)(1-2p)=1$. The coefficients of $\Xi_a$ have the signs $(+,+,\pm,-,-)$, because $0<n\le1$ gives $n^3-9n^2+6n+6\ge6-3n>0$ and $n^3-6n^2+2n+12\ge12-4n>0$. The rest is as for Proposition~\ref{prop:roots}.
\end{proof}

\begin{theorem}[Unique speed-driven orientation switch]\label{thm:cross}
There is exactly one speed $\hat v>1$ at which the best diagonal pair and the best axis pair have equal cost. The global optimum is diagonal for $1<v<\hat v$ and axis aligned for $v>\hat v$, and both are optimal at $\hat v$. There the optimal terminal locations jump. The optimal cost $F^*(v)$, the minimum of $F$ over terminal pairs, is continuous at $\hat v$, but its slope changes there.
\end{theorem}
\begin{proof}
Let $S_d^*(q)$ and $S_a^*(q)$ be the maximal savings of the two layouts. The roots in Propositions~\ref{prop:roots} and~\ref{prop:axisroots} are simple. By the implicit function theorem, the optimal offsets $a$ and $p$ are therefore smooth in $q$ on $[0,1)$, and so are $S_d^*$ and $S_a^*$. Appendix~\ref{app:cross} proves that $S_a^*/S_d^*$ is strictly decreasing in $q$. At $q=0$ this ratio exceeds one, because an axis pair then saves more than the diagonal pair with the same displacement $\delta_1$. To see this, note that $m=n=1$ at $q=0$, so these two pairs have $p=a=(1-\delta_1)/2$, $D_d=\delta_1$, and $D_a=\delta_1/2$, and \eqref{eq:Sd} and \eqref{eq:Sa} give $S_a-S_d=\delta_1^3(30-25\delta_1+\delta_1^2)/240>0$. As $q\uparrow1$ the ratio tends to zero, because $D_a\le(1-q)/2$ makes the axis saving vanish while the diagonal saving stays positive. The ratio therefore equals one at exactly one $q\in(0,1)$, with $S_d^*>S_a^*$ above it and $S_a^*>S_d^*$ below it, and by Theorem~\ref{thm:orient} the layout with the larger saving is globally optimal. At the tie $S_d^*=S_a^*>0$, so $F^*$ is continuous there. By \eqref{eq:logder}, $dS_d^*/dq>dS_a^*/dq$ at the tie, so the slope of $F^*$ changes there. The two optimal pairs differ, so the optimal locations jump.
\end{proof}
Because the ratio is monotone and $S_d^*-S_a^*$ is positive at $v=2.1830$ and negative at $v=2.1831$, the threshold satisfies $2.1830<\hat v<2.1831$. A diagonal link shortens a trip in both coordinates and so saves more per rider, but only trips that run its way in both coordinates can use it. An axis link shortens a trip in one coordinate only, but trips can use it whichever way they run in the other coordinate.

\section{Service implications: coverage, ridership, and the value of speed}\label{sec:service}
\subsection{Coverage and ridership}
\paragraph{Coverage.}
An origin is \emph{served} if a positive-area set of destinations offers a strict saving, and \emph{coverage} is the fraction of origins that are served. Put $\theta(X)=d(X,Q)-d(X,P)$, which measures how much closer the origin is to $P$ than to $Q$. By \eqref{eq:access}, no destination benefits when $|\theta(X)|\le L$. If $\theta(X)>L$, the destination $Q$ and its neighborhood benefit through $P$. The reverse direction applies when $\theta(X)<-L$. Thus an origin is served exactly when $|\theta(X)|>L$.

\begin{proposition}[Coverage of the two layouts]\label{prop:coverage}
The diagonal pair ${P=(a,a)}$, ${Q=(1-a,1-a)}$ and the axis pair ${P=(p,1/2)}$, ${Q=(1-p,1/2)}$, with offsets $0\le a,p<1/2$, serve these fractions of origins, respectively:
\begin{equation}\label{eq:coverage}
 2a^2+4aD_d+D_d^2\quad\text{with }D_d=m(1-2a)\text{ as in Section~\ref{sec:square}},\qquad\text{and}\qquad 1-(1-2p)/v\quad(v>1).
\end{equation}
\end{proposition}
\begin{proof}
For an axis pair, $\theta(X)=1-2X_1$ when $X_1$ lies between $p$ and $1-p$, and $|\theta(X)|=1-2p$ otherwise. The set $|\theta(X)|\le L$ is therefore a vertical strip of width $L=(1-2p)/v$. For a diagonal pair, coordinate $k$ contributes $|X_k-(1-a)|-|X_k-a|=(1-2a)-2\min\{[X_k-a]_+,1-2a\}$ to $\theta(X)$. Since $L=2(1-2a)-2D_d$ and $D_d\le1-2a$, we have $\theta(X)>L$ exactly when $[X_1-a]_++[X_2-a]_+<D_d$. This set consists of the square where both coordinates are at most $a$, two $a\times D_d$ rectangles, and a right triangle with legs $D_d$, so its area is $a^2+2aD_d+D_d^2/2$. Reflection through the center gives the set $\theta(X)<-L$, of equal area.
\end{proof}

\paragraph{Ridership.}
For fixed terminals let $L_0=d(X,Y)-\min\{d(X,P)+d(Q,Y),\,d(X,Q)+d(P,Y)\}$ be a trip's break-even link time, the time it would save on an instantaneous link. The trip strictly prefers the link exactly when $L_0>L$, and its travel time is $d(X,Y)-\pos{L_0-L}$. A marginal increase in $L$ therefore raises the travel time of exactly the trips with $L_0>L$, at unit rate. \emph{Ridership} is the fraction of these trips, $r=\Prob(L_0>L)$. With the terminals fixed, write $F(L)$ for the cost as a function of the link time. Dominated convergence gives its right and left derivatives,
\begin{equation}\label{eq:ridership}
 \lim_{\varepsilon\downarrow0}\frac{F(L+\varepsilon)-F(L)}{\varepsilon}=\Prob(L_0>L),\qquad
 \lim_{\varepsilon\downarrow0}\frac{F(L)-F(L-\varepsilon)}{\varepsilon}=\Prob(L_0\ge L)\quad(L>0).
\end{equation}

\begin{proposition}[Ridership of the two layouts]\label{prop:ridership}
For $q\in[0,1)$, the two pairs of Proposition~\ref{prop:coverage} carry these fractions of trips, respectively, with $D_d$ and $D_a$ as in Section~\ref{sec:square}:
\begin{equation}\label{eq:Pexplicit}
\begin{aligned}
 r_d&=\frac{D_d^4}{12}+\frac{4aD_d^3}{3}+6a^2D_d^2+8a^3D_d+2a^4,\\
 r_a&=p^2+2p(1+2p)D_a+\frac{1+8p-8p^2}{2}D_a^2+\frac{2(1-4p)}{3}D_a^3-\frac{D_a^4}{3}.
\end{aligned}
\end{equation}
\end{proposition}
\begin{proof}
Fix the terminals, that is, the offset $a$ of the diagonal pair or $p$ of the axis pair, and vary the speed. The link time $L=q\ell$ then varies, and by Theorem~\ref{thm:master} the cost is $2/3-S$ at every speed, where $S$ is \eqref{eq:Sd} or \eqref{eq:Sa} as a function $S_d(a,D)$ or $S_a(p,D)$ of the offset and of $D$. Since $dD/dL=-1/2$, the cost rises with $L$ at the rate $\frac12\,\partial S/\partial D$, which is the ridership by \eqref{eq:ridership}. Differentiating \eqref{eq:Sd} and \eqref{eq:Sa} at $D=D_d$ and $D=D_a$ gives \eqref{eq:Pexplicit}.
\end{proof}

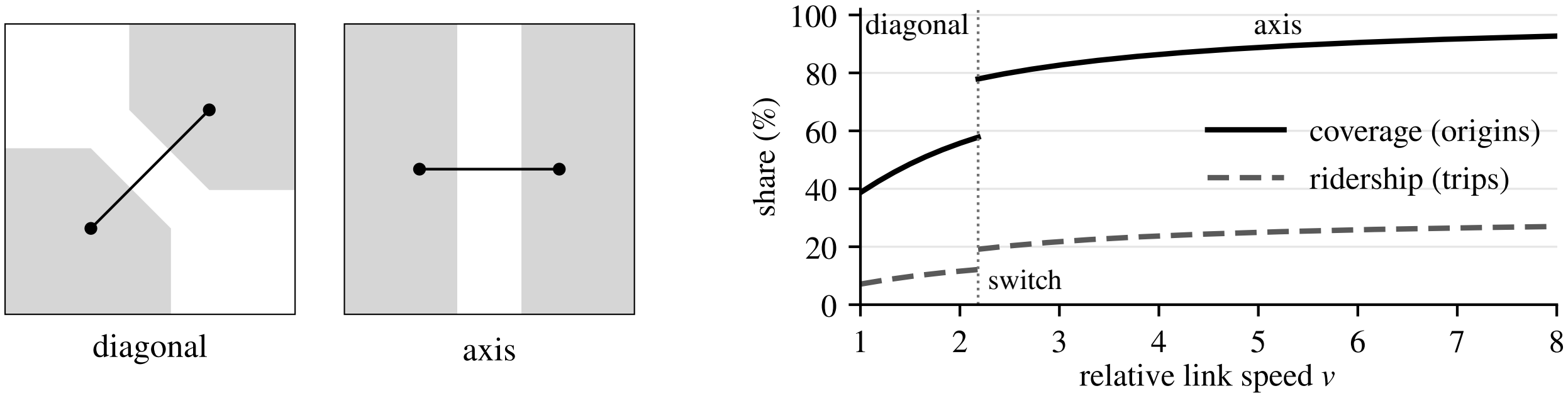
\begin{figure}[ht]
\setstretch{1.5}
\centering
\raisebox{3mm}{\begin{tikzpicture}[scale=3.05,font=\small]
\fill[black!16] (0,0)--(.5719,0)--(.5719,.2959)--(.2959,.5719)--(0,.5719)--cycle;
\fill[black!16] (1,1)--(.4281,1)--(.4281,.7041)--(.7041,.4281)--(1,.4281)--cycle;
\draw (0,0) rectangle (1,1);
\draw[thick] (.2959,.2959)--(.7041,.7041);
\fill (.2959,.2959) circle (.022);\fill (.7041,.7041) circle (.022);
\node at (.5,-.12) {diagonal};
\begin{scope}[xshift=1.17cm]
\fill[black!16] (0,0) rectangle (.3896,1);
\fill[black!16] (.6104,0) rectangle (1,1);
\draw (0,0) rectangle (1,1);
\draw[thick] (.2590,.5)--(.7410,.5);
\fill (.2590,.5) circle (.022);\fill (.7410,.5) circle (.022);
\node at (.5,-.12) {axis};
\end{scope}
\end{tikzpicture}}\hfill
\includegraphics{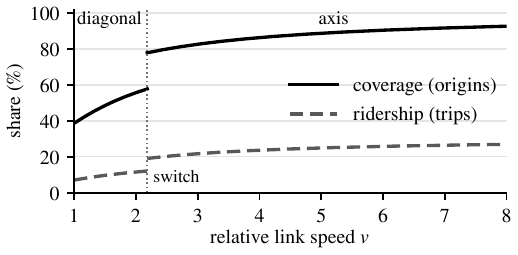}
\caption{Left: the two layouts that tie at $\hat v\approx2.18$, with served origins shaded. Right: coverage and ridership of the optimal design at each speed. Both jump at $\hat v$, whereas expected travel time is continuous.}\label{fig:service-profile}
\end{figure}

\paragraph{Equal cost, different service.}
Figure~\ref{fig:service-profile} shows the served origins of the two layouts that tie at $\hat v$ (left) and the coverage and ridership of the optimal design at each speed (right). The two corner squares of side $a$ away from a diagonal link are unserved at every speed, their points being equally far from both terminals. An axis link leaves unserved a strip whose width equals the link time $L$, so the strip vanishes as speed grows.

At $\hat v$, both layouts have expected travel time about $0.6339$, yet they serve different travelers. The diagonal layout covers $57.8\%$ of origins and carries $12.13\%$ of trips, and the axis layout $77.9\%$ and $19.12\%$. Equal total savings thus mean larger savings for fewer riders on the diagonal, about $0.270$ per rider against $0.172$. At every speed, Step~1 of Appendix~\ref{app:cross} gives a mean saving per rider of $\frac25(1-q/\sqrt2)$ for the optimal diagonal pair and at most $\frac13(1-q)$ for the optimal axis pair. Away from the tie the layouts stay close in travel time but not in service. At $v=2$, the best axis pair costs only $0.14\%$ more travel time than the optimal diagonal pair, yet it serves $76.5\%$ of origins instead of $55.7\%$ and $18.3\%$ of trips instead of $11.6\%$. A planner who also values coverage or ridership can adopt the axis layout below the switching speed at little cost in travel time.

\subsection{The value of speed and the choice of technology}
\paragraph{The marginal value of speed.}
At fixed terminals, a marginal increase in speed lowers $L=\ell/v$ at rate $\ell/v^2$, and by \eqref{eq:ridership} each unit of $L$ adds the ridership $r$ to expected travel time. Moving terminals that are already optimal has no first-order effect (envelope theorem). Ridership times link length therefore determines the marginal value of speed: with $r(v)$ and $\ell(v)$ those of the optimal design,
\begin{equation}\label{eq:value-speed}
 \frac{dF^*(v)}{dv}=-\frac{r(v)\ell(v)}{v^2}
 \qquad(v>1,\ v\ne \hat v).
\end{equation}
The gain from additional speed is thus proportional to the share of trips that ride the link.

\paragraph{Returns to speed diminish.}
As $v\to\infty$, the optimal cost $F^*(v)$ falls to about $0.5714$, so a single link saves at most $(2/3-0.5714)/(2/3)\approx14.3\%$ of ground travel time, and speeds five and ten already save $9.4\%$ and $11.6\%$. Although coverage then tends to one, ridership tends to only $30.7\%$, because access and egress still make ground travel preferable for most trips.

\paragraph{The choice of technology.}
Suppose now that speed has a price and is chosen together with the terminals.

\begin{corollary}[Joint choice of speed and alignment]\label{cor:joint}
Let the planner minimize $F^*(v)+c(v)$ over $v\ge1$, where $v=1$ means ground speed, $F^*(1)$ is the limit of $F^*(v)$ as $v\downarrow1$, and $c$ is a speed cost per trip in normalized time. If $c$ is differentiable at $\hat v$, that speed is never optimal. Let $c(v)=\gamma(v-1)$, where $\gamma>0$ is the price of a unit of speed, so that ground speed costs nothing. Then there is a threshold price $\hat\gamma$: if $\gamma>\hat\gamma$, every optimal design is a diagonal link; if $\gamma<\hat\gamma$, every optimal design is an axis-aligned link; at $\gamma=\hat\gamma$ both are optimal. There are speeds $v_-<\hat v<v_+$ such that every optimal diagonal link has speed at most $v_-$ and every optimal axis-aligned link has speed at least $v_+$, so no speed strictly between $v_-$ and $v_+$ is optimal at any price.
\end{corollary}
\begin{proof}
The first step treats a differentiable cost, and the other three the linear cost ${c(v)=\gamma(v-1)}$.

\emph{The switching speed is never optimal:} The derivation of \eqref{eq:value-speed} applies to each layout separately, so its maximal saving satisfies $dS^*/dq=-r\ell$, with its own ridership and link length. At the tie $dS_d^*/dq>dS_a^*/dq$ (proof of Theorem~\ref{thm:cross}), so $r_a\ell_a>r_d\ell_d$. Since the axis layout takes over at $\hat v$, the right derivative of $F^*$ there, $-r_a\ell_a/v^2$, is smaller than its left derivative, $-r_d\ell_d/v^2$. The same holds for $F^*+c$, whereas a minimum requires a nonpositive left derivative and a nonnegative right derivative.

\emph{Optimal speeds fall as the price rises:} Let $J_\gamma(v)=F^*(v)+\gamma(v-1)$. Its minimizers exist and form a compact set, because $J_\gamma$ is continuous, $J_\gamma(1)\le2/3$, and $J_\gamma(v)>2/3$ for $v>1+2/(3\gamma)$. If $\gamma_1<\gamma_2$ and $v_i$ minimizes $J_{\gamma_i}$, then $J_{\gamma_1}(v_1)\le J_{\gamma_1}(v_2)$ and $J_{\gamma_2}(v_2)\le J_{\gamma_2}(v_1)$. Adding these inequalities cancels the travel-time terms and leaves $(\gamma_2-\gamma_1)(v_2-v_1)\le0$, so $v_2\le v_1$.

\emph{Both alignments occur:} Fix $\bar v>\hat v$. By \eqref{eq:value-speed}, $F^*$ is nonincreasing and $F^*(\bar v)<F^*(\hat v)$. Every $v\le\hat v$ has $J_\gamma(v)\ge F^*(\hat v)$, which exceeds $J_\gamma(\bar v)=F^*(\bar v)+\gamma(\bar v-1)$ once $\gamma(\bar v-1)<F^*(\hat v)-F^*(\bar v)$. At such low prices every optimal speed therefore exceeds $\hat v$. Conversely, every $v\ge\hat v$ has $J_\gamma(v)\ge\gamma(\hat v-1)$, which exceeds $2/3\ge J_\gamma(1)$ once $\gamma(\hat v-1)>2/3$. At such high prices every optimal speed is therefore below $\hat v$.

\emph{Threshold and gap:} Let $\hat\gamma$ be the supremum of the prices at which some optimal speed exceeds $\hat v$. By the last step, $0<\hat\gamma<\infty$. Because optimal speeds fall as the price rises, every optimal speed exceeds $\hat v$ when $\gamma<\hat\gamma$ and none does when $\gamma>\hat\gamma$, and by the first step none equals $\hat v$. By Theorem~\ref{thm:cross}, optimal designs are therefore axis aligned below $\hat\gamma$ and diagonal above it. Since $J_\gamma(v)$ is continuous in $(\gamma,v)$, limits of optimal speeds as $\gamma\uparrow\hat\gamma$ and as $\gamma\downarrow\hat\gamma$ are optimal at $\hat\gamma$, on the two sides of $\hat v$, so both layouts are optimal there. Let $v_+$ be the smallest optimal speed above $\hat v$ at $\hat\gamma$ and $v_-$ the largest one below. No optimal speed lies between them. By the second step, optimal speeds are at most $v_-$ at higher prices and at least $v_+$ at lower ones.
\end{proof}
Some technologies are never worth building, because no optimal design uses a speed in a band around $\hat v$, whatever the price of speed. The choice of technology is lumpy. As speed becomes cheaper, the planner does not accelerate a diagonal link gradually but replaces it by a markedly faster axis link at other sites, which serves more origins and carries more trips. Speed and location should therefore be decided together.

\section{Extension: a rectangular district}\label{sec:rectangle}
Many districts, such as waterfronts and airport corridors, are elongated. Let the district be $\Omega_\rho=[0,\rho]\times[0,1]$ with $\rho>1$, so that expected ground distance is $(1+\rho)/3$, and write $F_\rho$ for expected travel time with the link. Lemma~\ref{lem:window}, scaled to each side, preserves the centering proof, so every optimum is centered. Write its displacement as $(x,y)\in[0,\rho]\times[0,1]$, where $x$ is the displacement along the long side and $y$ the one along the short side. As on the square, $\delta_1=\max(x,y)$ is the larger of the two displacements and $\delta_2=\min(x,y)$ the smaller, but the sides that carry them now differ in length. Let $s_1$ and $s_2$ be the lengths of the sides that carry $\delta_1$ and $\delta_2$, so that $(s_1,s_2)=(\rho,1)$ if $x\ge y$ and $(1,\rho)$ if $y>x$. Let $A=(s_1-\delta_1)/2$, $b=(s_2-\delta_2)/2$, $D=(\delta_1+\delta_2-L)/2$, and $E=[\delta_1-\delta_2-L]_+/2$. Let $G$ be the forward saving $[g_{PQ}]_+$ integrated over all origins and destinations in $\Omega_\rho$. Because $\Omega_\rho$ has area $\rho$, the expected forward saving is $G/\rho^2$. A rectangle allows one case that the square does not. A placement is \emph{deep} if $\delta_1>s_2+L$, which requires the larger displacement to lie along the long side, and we put $E'=[\delta_1-s_2-L]_+/2$.

\begin{proposition}[Exact objective on a rectangle]\label{prop:rectformula}
For every centered terminal placement on $\Omega_\rho$ and every link time $0\le L\le\delta_1+\delta_2$,
\begin{equation}\label{eq:rectformula}
 F_\rho=\frac{1+\rho}{3}-\frac{2G}{\rho^2},\qquad
 G=R(A,b,D)+C_{A,b,s_2}(E)+\frac{4A^2E'^3}{3}+\frac{2AE'^4}{3}+\frac{E'^5}{15}.
\end{equation}
\end{proposition}
\begin{proof}
We repeat the proof of Theorem~\ref{thm:master} with areas instead of probabilities. On a side of length $s$, the trips of one coordinate fill a square of area $s^2$, and centered terminals separated by $\delta$ leave the gap $(s-\delta)/2$ at each end. The proof of Lemma~\ref{lem:window} then gives the mass at $\delta$ and the central density of \eqref{eq:law} with this gap in place of $(1-\delta)/2$, the mass $s^2/2-(s-\delta)^2/4$ at $-\delta$, and the lower density $s+z$ on $-s<z<-\delta$. The function $\Phi$ and the regular expression are therefore unchanged, and in \eqref{eq:master-integral} the mass at $-\delta_2$ becomes $s_2^2/2-b^2$ and the lower density becomes $s_2+z_2$. This gives $G=R(A,b,D)+C_{A,b,s_2}(E)$ as long as the threshold $L-\delta_1$ stays in the range of the lower density, that is, $\delta_1-L\le s_2$. For a deep placement the lower density ends at $z_2=-s_2$, above the threshold, whereas the formula continues it as $s_2+z_2$ on $L-\delta_1<z_2<-s_2$, an interval of length $2E'$ on which the true density is zero. Removing this part adds $\int_{L-\delta_1}^{-s_2}(-s_2-z_2)\Phi(z_2)\dd z_2=4A^2E'^3/3+2AE'^4/3+E'^5/15$. As in \eqref{eq:saving}, the two directions together save $2G/\rho^2$ per trip.
\end{proof}

\begin{theorem}[Optimal layouts on a rectangle]\label{thm:rectglobal}
Every optimum is centered. For $1<v<\infty$, an optimal displacement lies on the long axis, $y=0$, or is an interior stationary point of the regular regime; all such points have $0<x-y<\rho-1$.
\end{theorem}
\begin{proof}
Centering was shown above. An optimum with $y=0$ is on the long axis. By Appendix~\ref{app:exclude}, one with $y>0$ has $y\le x$ (Part~1), is not irregular (Parts~2 and~3), and lies neither on the interface $L=x-y$ nor on the faces $x=\rho$ and $y=1$ (Part~3). It is therefore an interior point of the regular regime, where $G=R(A,b,D)$ is differentiable, and hence stationary. Here $D=(x+y-q\ell)/2$ with $\ell=\sqrt{x^2+y^2}$, and $A=(\rho-x)/2$ and $b=(1-y)/2$ if $x\ge y$. If $y>x$, the two exchange values, which does not change $R$, because $R$ is symmetric in its first two arguments. Regard $R$ as a function of $A+b$, $\alpha=2Ab$, and $D$, as in the proof of Theorem~\ref{thm:orient}. The $x$- and $y$-derivatives of $A+b$ are equal, those of $\alpha$ are $-b$ and $-A$, and those of $D$ are $(1-qx/\ell)/2$ and $(1-qy/\ell)/2$. Subtracting the two first-order conditions therefore gives ${q(x-y)R_D=\ell\{\rho-1-(x-y)\}R_\alpha}$, where $R_\alpha,R_D>0$ because $D>0$ when $x,y>0$. If $x\le y$, the left side would be nonpositive and the right side positive. Hence $x>y$, so the right side is positive as well and $x-y<\rho-1$.
\end{proof}
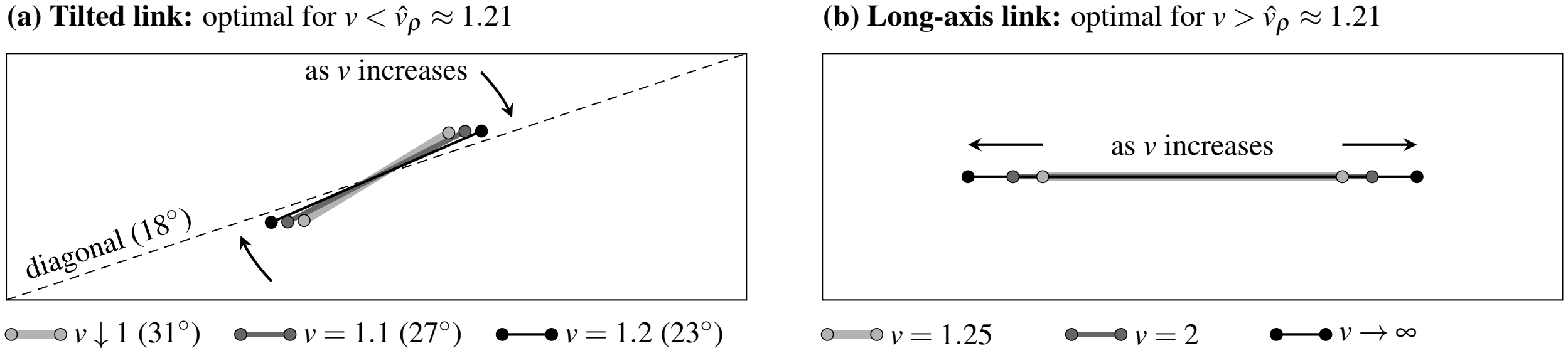
\begin{figure}[t]
\setstretch{1.5}
\centering
\begin{tikzpicture}[x=2.6cm,y=2.6cm,font=\footnotesize,>=stealth]
\tikzset{la/.style={line width=2.6pt,black!30},lb/.style={line width=1.6pt,black!60},lc/.style={line width=0.8pt,black},
 ta/.style={circle,draw=black,line width=.3pt,fill=black!30,inner sep=0pt,minimum size=3.6pt},tb/.style={ta,fill=black!60},tc/.style={ta,fill=black},
 turn/.style={->,line width=0.8pt}}
\node[anchor=south west,inner sep=0pt] at (0,1.07) {\textbf{(a) Tilted link:} optimal for $v<\hat v_\rho\approx1.21$};
\draw[line width=.4pt,densely dashed] (0,0)--(3,1);
\node[rotate=18.43,anchor=south,inner sep=1.5pt] at (0.44,0.1467) {diagonal ($18^\circ$)};
\draw[la] (1.2068,0.3221)--(1.7932,0.6779);
\draw[lb] (1.1408,0.3165)--(1.8592,0.6835);
\draw[lc] (1.074,0.3148)--(1.926,0.6852);
\node[ta] at (1.2068,0.3221) {};\node[ta] at (1.7932,0.6779) {};
\node[tb] at (1.1408,0.3165) {};\node[tb] at (1.8592,0.6835) {};
\node[tc] at (1.074,0.3148) {};\node[tc] at (1.926,0.6852) {};
\draw[turn] (1.0757,0.0758) arc[start angle=225.00,end angle=203.50,radius=0.6];
\draw[turn] (1.9243,0.9242) arc[start angle=45.00,end angle=23.50,radius=0.6];
\node[anchor=east,inner sep=0pt] at (1.8743,0.9342) {as $v$ increases};
\draw (0,0) rectangle (3,1);
\draw[la] (0.02,-.14)--(0.22,-.14);\node[ta] at (0.02,-.14) {};\node[ta] at (0.22,-.14) {};\node[anchor=west,inner sep=0pt] at (0.27,-.14) {$v\downarrow1$ ($31^\circ$)};
\draw[lb] (0.95,-.14)--(1.15,-.14);\node[tb] at (0.95,-.14) {};\node[tb] at (1.15,-.14) {};\node[anchor=west,inner sep=0pt] at (1.2,-.14) {$v=1.1$ ($27^\circ$)};
\draw[lc] (2.01,-.14)--(2.21,-.14);\node[tc] at (2.01,-.14) {};\node[tc] at (2.21,-.14) {};\node[anchor=west,inner sep=0pt] at (2.26,-.14) {$v=1.2$ ($23^\circ$)};
\begin{scope}[xshift=8.6cm]
\node[anchor=south west,inner sep=0pt] at (0,1.07) {\textbf{(b) Long-axis link:} optimal for $v>\hat v_\rho\approx1.21$};
\draw[la] (0.8929,0.5)--(2.1071,0.5);
\draw[lb] (0.7724,0.5)--(2.2276,0.5);
\draw[lc] (0.5903,0.5)--(2.4097,0.5);
\node[ta] at (0.8929,0.5) {};\node[ta] at (2.1071,0.5) {};
\node[tb] at (0.7724,0.5) {};\node[tb] at (2.2276,0.5) {};
\node[tc] at (0.5903,0.5) {};\node[tc] at (2.4097,0.5) {};
\draw[turn] (2.1071,0.63)--(2.4097,0.63);\draw[turn] (0.8929,0.63)--(0.5903,0.63);
\node[inner sep=0pt] at (1.5,0.63) {as $v$ increases};
\draw (0,0) rectangle (3,1);
\draw[la] (0.02,-.14)--(0.22,-.14);\node[ta] at (0.02,-.14) {};\node[ta] at (0.22,-.14) {};\node[anchor=west,inner sep=0pt] at (0.27,-.14) {$v=1.25$};
\draw[lb] (1.01,-.14)--(1.21,-.14);\node[tb] at (1.01,-.14) {};\node[tb] at (1.21,-.14) {};\node[anchor=west,inner sep=0pt] at (1.26,-.14) {$v=2$};
\draw[lc] (1.84,-.14)--(2.04,-.14);\node[tc] at (1.84,-.14) {};\node[tc] at (2.04,-.14) {};\node[anchor=west,inner sep=0pt] at (2.09,-.14) {$v\to\infty$};
\end{scope}
\end{tikzpicture}
\caption{Optimal link on the rectangle with $\rho=3$ at six speeds, found by minimizing \eqref{eq:rectformula} numerically. The tilted link and the long-axis link tie at the speed $\hat v_\rho\approx1.21$. Angles are those with the long side.}\label{fig:rectlinks}
\end{figure}
Only two families of layouts can thus be optimal on a rectangle. The regular candidate takes the place of the square's diagonal pair. Its displacement along the long side is the larger one, but by less than the difference $\rho-1$ of the side lengths, so neither a link at $45$ degrees nor one along the short axis is optimal. The planner compares this candidate with the best long-axis link, found by a one-dimensional search. Figure~\ref{fig:rectlinks} shows the optimal link for $\rho=3$ at six speeds. As $v$ increases, the tilted link lengthens and turns toward the long side, though not as far as the diagonal. On nearly square rectangles ($\rho\le1.31$ in our computations), it turns past the diagonal before the long axis takes over, but by at most $0.4$ degrees. Our computations ($1.01\le\rho\le10$) suggest that each rectangle has only one switching speed, but we cannot prove this analytically.

\section{Conclusions and implementation scope}\label{sec:conclusions}
We have studied the placement of an optional fast connection whose users must reach its terminals and may instead travel on a slower ground network. For uniform demand on a square, the two-terminal problem has a complete global analytical solution. Every optimal terminal pair is symmetric about the center, only diagonal and axis-aligned layouts can be optimal, and the orientation changes from diagonal to axis aligned at a single relative speed. Better technology can thus move an optimal facility, not merely speed travel through it. A planner who chooses technology and terminals together, with a differentiable cost of speed, never builds at the switching speed, and with a cost linear in speed skips a band of speeds around it. The layouts tied at the threshold serve different shares of the region and trips, so evaluations should report the distribution of benefits alongside average travel time. Two practical lessons follow. First, near the switching speed, coverage or ridership can decide between the two layouts at almost no cost in travel time. Second, returns to speed diminish, and one link saves at most about $14\%$ of travel time. In a rectangular district, the optimal connection runs along the long axis or is tilted toward the long side.

Our model is an analytical benchmark. A plan for a real service would have to account for observed demand, boarding and waiting times, the sites available for terminals, costs, and capacity, and the planner might also weigh reliability or the distribution of benefits. We have not studied networks with more than two terminals, in which several links can serve the same trip, and leave them for future research.

\paragraph{Code and data availability.}
No empirical data are used, and every numbered result is proved in the paper. The reported numbers and Figure~\ref{fig:service-profile} are computed from the closed forms, and Figure~\ref{fig:rectlinks} from a numerical minimization of \eqref{eq:rectformula}. The code that reproduces them is provided as a supplementary file.

\appendix
\section{Excluding other orientations}\label{app:exclude}
This appendix proves the exclusions used in Theorems~\ref{thm:orient} and \ref{thm:rectglobal}, on the rectangle $\Omega_\rho$ with $\rho\ge1$. The square is the case $\rho=1$. Let $(x,y)$ be the displacements of a centered pair along the long and the short side, and let $G$ be the saving in \eqref{eq:rectformula}. Subscripts denote partial derivatives with $L$ held fixed, whereas $dG/dx$ and $dG/dy$ include the dependence of $L=q\ell$ on $\ell=\sqrt{x^2+y^2}$.

\paragraph{Part 1: the larger displacement belongs on the long side.}
Let $\rho>1$, and suppose a design assigns $\delta_1>\delta_2$ to the short side, so $\delta_1\le1$. Exchanging the two displacements changes neither $L$ nor $D$. The two gaps are $(1-\delta_1)/2$ and $(\rho-\delta_2)/2$ before the exchange and $(\rho-\delta_1)/2$ and $(1-\delta_2)/2$ after it, so $A+b$ is unchanged and $\alpha=2Ab$ rises by $(\rho-1)(\delta_1-\delta_2)/2$. A regular placement therefore gains, because $G=R$ depends on $A$ and $b$ only through $A+b$ and $\alpha$ and has $R_\alpha>0$ (proof of Theorem~\ref{thm:orient}). For an irregular placement, $E=(\delta_1-\delta_2-L)/2>0$, and neither assignment is deep, since $\delta_1\le1$. Write $B=(1+\rho-2\delta_1)/2$ and $H=(1-\delta_1)(\rho-\delta_1)/4$, which are nonnegative. Inserting $(A,b,s)=((1-\delta_1)/2,(\rho-\delta_2)/2,\rho)$ for the short assignment and $((\rho-\delta_1)/2,(1-\delta_2)/2,1)$ for the long one into \eqref{eq:R}--\eqref{eq:C}, with $D=E+\delta_2$ and $\delta_1=2E+\delta_2+L$, gives
\begin{align*}
 G_{\rm long}-G_{\rm short}=\frac{\rho-1}{12}\bigl[&26E^4+(28L+36B+52\delta_2)E^3+(6L^2+36LB+42L\delta_2+48H+54B\delta_2+30\delta_2^2)E^2\\*
 &+(6L^2B+6L^2\delta_2+24LH+36LB\delta_2+18L\delta_2^2+48H\delta_2+18B\delta_2^2+4\delta_2^3)E\\*
 &+3L^2B\delta_2+3L^2\delta_2^2+12LH\delta_2+6LB\delta_2^2+2L\delta_2^3\bigr].
\end{align*}
All terms in brackets are nonnegative and the first is positive, so the difference is positive. At $\rho=1$ the exchange is a symmetry of the square. Hence we may take $y\le x$ in both cases.

\paragraph{Part 2: no tilted irregular placement is stationary.}
Now $\delta_1=x$ and $\delta_2=y$, so $A=(\rho-x)/2$ and $b=(1-y)/2$. Put $h=x-L$. Then $D=(h+y)/2$, $E=[h-y]_+/2$, and $E'=[h-1]_+/2$, so a tilted irregular placement has $0<y<h$, and it is deep when $h>1$. Differentiating \eqref{eq:rectformula} at fixed $L$ gives
\begin{equation}\label{eq:NG}
 24G_y=yN,\qquad
 N=-24A^2-24Ah+4Ay^2+12A-3h^2+hy^2+3h+y^3-3y^2,
\end{equation}
also for deep placements, whose additional terms do not depend on $y$. Since $\partial L/\partial y=Ly/\ell^2$, we have $dG/dy=yN/24+(Ly/\ell^2)G_L$, where $G_L\le0$ because more link time cannot raise the saving. Hence reducing $y$ improves the placement whenever $N<0$.

If $h\ge1$, then $N=-24A^2+4A(y^2-6h+3)+(h-1)(y^2-3h)+y^2(y-2)<0$, because no term is positive and the last is negative. If $h<1$, let $\eta=1-h-2A$, which equals $L-(\rho-1)\le L$. In terms of $\eta$,
\[
 N=-6\eta^2+(6-2y^2)\eta-3h(1-h)-y^2(1+h-y),
\]
so that $N<0$ whenever $\eta\le0$, in particular if $q=0$.

It remains to exclude a stationary point with $h<1$ and $q>0$. By \eqref{eq:rectformula}, $96G_x=K$ and $96G_L=-M$, where
\begin{align*}
 T&=h^2(6+4h-h^2)+(6-12h)y^2+y^4,\\
 K&=48A^2(1+2h-h^2)-T,\\
 M&=48A^2(1+2h-h^2)+16Ah(3+3h-h^2)-48Ay^2+T.
\end{align*}
Here $T>0$, because its terms are nonnegative if $h\le1/2$, and otherwise $y<h$ gives $T>h^2(12-8h-h^2)>0$. Also $16Ah(3+3h-h^2)-48Ay^2\ge16Ah(3-h^2)\ge0$, so $M>0$. The stationarity conditions $dG/dx=dG/dy=0$ read $\ell^2K=LxM$ and $4\ell^2N=LM$. Hence $N>0$, so that $\eta>0$, and $K=4xN$. Put $V=M-K-4hN=M-4(x+h)N$. Since $L=x-h$ and $LM=4(x^2+y^2)N$, we have ${LV=LM-4(x^2-h^2)N=4(y^2+h^2)N}$, so that $V>0$. Expansion with $\eta=1-h-2A$ gives the identity
\begin{align*}
 \eta V-4(y^2+h^2)N={}&4(h^2+y^2)[12A^2h+2A(3h^2+y^2)+y^2(2h-y)]\\*
 &+2\eta[96A^3h+72A^2(2h^2-y^2)+4Ah(15h^2-11y^2)+3h^2(h^2-y^2)+y^2(h-y)^2+6y^4]\\*
 &+8\eta^2[12A^2h+6A(2h^2-y^2)+h^3]+12\eta^3y^2.
\end{align*}
Every bracket is nonnegative and the first is positive, because $0<y<h<1$, $A\ge0$, and $\eta>0$. Therefore $\eta V>4(y^2+h^2)N=LV$, that is, $\eta>L$, which contradicts $\eta\le L$.

\paragraph{Part 3: no optimum with $y>0$ lies on the interface or on an outer face.}
Consider first a point on the interface $L=x-y$ with $A>0$. The correction $C_{A,b,1}(E)$ has the linear term $2A^2(\frac12-b^2)E$, which is positive because $b\le1/2$, and $E=[x-y-L]_+/2$ grows at rate $(1-qx/\ell)/2>0$ as $x$ increases through the interface. The derivative of $G$ in $x$ therefore jumps upward there, and a point with such a kink cannot maximize $G$. If $A=0$, the point lies on the face $x=\rho$. On the face $y=1$ we have $b=0$, and on the regular side
\[
 2\,\frac{dG}{dy}=\Bigl(1-\frac{qy}{\ell}\Bigr)R_D-R_b\le R_D-R_b=-\frac{D^4}{12}-\frac{2AD^3}{3}-A^2D^2<0,
\]
so decreasing $y$ improves the design. On the irregular side $h>1$, and Part~2 gives the same conclusion. On the face $x=\rho$ with $y>0$ we have $A=0$. On the regular side, which includes the interface, the same computation with $A$ and $b$ exchanged gives ${2\,dG/dx\le R_D-R_A=-D^4/12-2bD^3/3-b^2D^2<0}$, so decreasing $x$ improves the design. On the irregular side with $h<1$, Part~2 gives $96G_x=K=-T<0$ and $G_L\le0$, so $dG/dx=G_x+(Lx/\ell^2)G_L<0$ and decreasing $x$ improves the design again. If $h\ge1$, decreasing $y$ improves it by Part~2.

\section{Proof of the unique orientation switch}\label{app:cross}
We prove that $S_a^*/S_d^*$ decreases strictly in $q$ on $[0,1)$, as the proof of Theorem~\ref{thm:cross} requires. Throughout, $S_d(a,D)$ and $S_a(p,D)$ denote \eqref{eq:Sd} and \eqref{eq:Sa} with $D$ free, as in Section~\ref{sec:service}, and $a$, $p$, $\lambda_d$, and $\lambda_a$ are the optimal offsets and the roots of Propositions~\ref{prop:roots} and~\ref{prop:axisroots}, so that $1-2a=1/(1+\lambda_d)$ and $1-2p=1/(1+\lambda_a)$.

\paragraph{Step 1: the rates at which the maximal savings shrink.}
The diagonal saving is homogeneous of degree five in $(a,D)$, whereas the axis saving is a sum of terms of degrees three, four, and five in $(p,D)$. Euler's theorem on homogeneous functions therefore gives $a\,\partial S_d/\partial a+D\,\partial S_d/\partial D=5S_d$ and $p\,\partial S_a/\partial p+D\,\partial S_a/\partial D=3S_a+\chi$, where $\chi$ counts the terms of degree four once and those of degree five twice:
\[
 \chi=4p^2D^2\Bigl(1-\frac{4D}{3}\Bigr)+\frac{8pD^3}{3}(1-D)+D^4\Bigl(\frac13-\frac{4D}{15}\Bigr)\ge0\qquad(D\le\tfrac12).
\]
At the optimal offsets, $D_d=m(1-2a)$ and $D_a=n(1-2p)/2$, and the first-order conditions read $\partial S_d/\partial a=2m\,\partial S_d/\partial D$ and $\partial S_a/\partial p=n\,\partial S_a/\partial D$. Substituting them into the two identities obtained from Euler's theorem, with $2ma+D_d=m$ and $np+D_a=n/2$, gives $m\,\partial S_d/\partial D=5S_d$ and $(n/2)\,\partial S_a/\partial D=3S_a+\chi$. With $r=\frac12\,\partial S/\partial D$ from Proposition~\ref{prop:ridership}, these two relations give the savings per rider used in Section~\ref{sec:service}, $S_d/r_d=2m/5$ and $S_a/r_a=nS_a/(3S_a+\chi)\le n/3$. By the envelope theorem, $q$ changes the maximal savings only through $D_d$ and $D_a$, at the rates $\partial D_d/\partial q=-(1-2a)/\sqrt2$ and $\partial D_a/\partial q=-(1-2p)/2$. With $\sqrt2\,m=\sqrt2-q$ and $n=1-q$, the relative rates are therefore
\begin{equation}\label{eq:scaling}
 -\frac{d\log S_d^*}{dq}=\frac{1-2a}{\sqrt2\,S_d}\,\frac{\partial S_d}{\partial D}=\frac{5(1-2a)}{\sqrt2-q},\qquad
 -\frac{d\log S_a^*}{dq}=\frac{1-2p}{2S_a}\,\frac{\partial S_a}{\partial D}=\frac{1-2p}{1-q}\left(3+\frac{\chi}{S_a}\right).
\end{equation}

\paragraph{Step 2: bounds on the optimal offsets.}
By Propositions~\ref{prop:roots} and~\ref{prop:axisroots}, each of $\Xi_d$ and $\Xi_a$ is negative below its positive root and positive above it. Direct substitution, with $m=1-q/\sqrt2$ and $n=1-q$, gives
\begin{align*}
 -1728\,\Xi_d\bigl(\tfrac56+\tfrac{3q}{2\sqrt2}\bigr)
 &=2943q^4+14292\sqrt2\,q^3+21060q^2-888\sqrt2\,q+476>0,\\*
 8748\,\Xi_a\bigl(\tfrac{11+25q}{18}\bigr)
 &=(1-q)(390773q^3+715485q^2+358395q+5011)>0,\\*
 \Xi_a\bigl(\tfrac12\bigr)
 &=-n^4-(1-n)\bigl(4n^2+\tfrac{49}2n+\tfrac{31}2\bigr)-\tfrac{19}4<0.
\end{align*}
In the first line, the last three terms form a positive quadratic because $(888\sqrt2)^2<4\cdot21060\cdot476$. Therefore $\lambda_d>5/6+3q/(2\sqrt2)$ and $1/2<\lambda_a<(11+25q)/18$. Since $1-2a=1/(1+\lambda_d)$, $1-2p=1/(1+\lambda_a)$, and $p=\lambda_a(1-2p)/2$, these bounds give
\begin{equation}\label{eq:optbounds}
 1-2a<\frac{12}{22+9\sqrt2\,q},\qquad
 1-2p>\frac{18}{29+25q},\qquad
 D_a=\frac{n(1-2p)}{2}<\frac13,\qquad \frac p{D_a}=\frac{\lambda_a}{n}>\frac12.
\end{equation}
The last two bounds use $\lambda_a>1/2$, which gives $1-2p<2/3$, and $n\le1$.

\paragraph{Step 3: a lower bound on $\chi/S_a$.}
For $0<D\le1/3$ and $p\ge D/2$, expansion gives
\begin{align*}
 5\chi-3DS_a
 &=D^2(14p^2+22pD/3+2D^2/3)-D^3(116p^2/3+64pD/3+7D^2/3)\\*
 &\quad+D^4(8p^2+4pD+2D^2/5)
 \ \ge\ D^2(10p^2+2pD-D^2)/9>0.
\end{align*}
The bound uses $D\le1/3$ in the negative term and drops the last positive term, and the final inequality uses $p\ge D/2$. By \eqref{eq:optbounds}, both conditions hold at the axis optimum, so $\chi/S_a\ge3D_a/5=3n(1-2p)/10$.

\paragraph{Step 4: the axis saving shrinks faster.}
By Step~3, the axis rate is at least ${(1-2p)\{3+3n(1-2p)/10\}/n}$, which increases in $1-2p$. Inserting the bounds on $1-2p$ and $1-2a$ from \eqref{eq:optbounds} into the two rates gives
\begin{equation}\label{eq:logder}
 \frac{d}{dq}\log\frac{S_d^*}{S_a^*}
 >\frac{108(77+58q)}{5(29+25q)^2(1-q)}
 -\frac{60}{(\sqrt2-q)(22+9\sqrt2\,q)}.
\end{equation}
Over the common denominator $5(29+25q)^2(1-q)(\sqrt2-q)(22+9\sqrt2\,q)$, which is positive, the right side has the numerator $12\,\Gamma(q)$, where $\Gamma(q)=(15246\sqrt2-21025)-(17997-11484\sqrt2)q+(18537-6237\sqrt2)q^2+(15625-4698\sqrt2)q^3$. Since $1.414<\sqrt2<1.415$, $\Gamma(q)>530-1760q+9700q^2$ for $q\ge0$, and this quadratic is positive because $1760^2<4\cdot530\cdot9700$. Hence the right side of \eqref{eq:logder} is positive, and $S_a^*/S_d^*$ decreases strictly.

\label{LastMainPage}
\setlength{\bibsep}{3pt}

\end{document}